\documentclass[AMA,STIX1COL]{WileyNJD-v2}

\usepackage{moreverb}

\newcommand\BibTeX{{\rmfamily B\kern-.05em \textsc{i\kern-.025em b}\kern-.08em
T\kern-.1667em\lower.7ex\hbox{E}\kern-.125emX}}

\articletype{reasearch article}%

\received{<day> <Month>, <year>}
\revised{<day> <Month>, <year>}
\accepted{<day> <Month>, <year>}

\begin{document}

\title{Stochastic Maximum Principle for a generalized Volterra Control System}

\author[1]{Yuhang Li}

\author[2]{Yuecai Han}

\authormark{Li and Han }

\address{\orgdiv{School of Mathematics}, \orgname{Jilin University}, \orgaddress{\state{Changchun 130012}, \country{China}}}

\corres{Yuecai Han, School of Mathematics, Jilin University,
Changchun 130012, China. \email{hanyc@jlu.edu.cn}}

\presentaddress{School of Mathematics, Jilin
University,  Changchun 130012,  China}

\abstract[Abstract]{In this paper, we consider the stochastic optimal control problem for a generalized Volterra control system. The corresponding state process  is a kind of  a generalized stochastic Volterra integral differential equations. We prove the existence and uniqueness of the solution of this type of equations. We obtain the stochastic maximum principle of the  optimal control system by introducing a kind of generalized anticipated backward stochastic differential equations. We prove the existence and uniqueness of the solution of this adjoint equation, which may be singular at some points. As an application,  the linear quadratic control problem is investigated  to illustrate the main results.}


\maketitle

\footnotetext{\textbf{Abbreviations:} ANA, anti-nuclear antibodies; APC, antigen-presenting cells; IRF, interferon regulatory factor}

\section{Introduction}\label{sec1}

\quad To better describe the real-world, the stochastic integral differential equations have been studied in many areas, such as in biological science, applied mathematics, physics, and other disciplines, etc \cite{bloom1980bounds,holmaaker1993global,forbes1997caluculating,du2014reproducing}.
  Mao and Riedle \cite{mao2006mean} study the stability of some types of stochastic volterra integral differential equations. Nesterenko\cite{nesterenko2014modified}  substantiate the application of a modified projection-iterative method to the solution of boundary value problems for weakly nonlinear integrodifferential equations with parameters. Dzhumabaev\cite{dzhumabaev2016one} establish the necessary and sufficient conditions for the well-posedness of linear boundary value
problems for Fredholm integro-differential equations . Zhang et al.\cite{zhang2020theoretical} investigate numerical analysis of the following generalized stochastic volterra integral differential equations
 \begin{align*}
d Y(t)=&f\left(Y(t), \int_0^t k_1(t, s) Y(s) d s, \int_0^t \sigma_1(t, s) Y(s) d w(s)\right) d t\\
&+g\left(Y(t), \int_0^t k_2(t, s) Y(s) d s, \int_0^t \sigma_2(t, s) Y(s) d w(s)\right) d w(t).
 \end{align*}
They prove the existence and uniqueness of the solution when $\|k_i\|$ and $\|\sigma_i\|$ are bounded.

Control problems for integral differential equations have also been studied. Kim\cite{kim1993control}  discuss a reachability problem for a second-order
integro-differential equation based
upon a new kind of unique continuation property. Mashayekhi et al.\cite{mashayekhi2013hybrid} give a new numerical method for solving the optimal control of a class of systems
described by integro-differential equations with quadratic performance index. Assanova et al.\cite{assanova2020numerical} present the  existence of optimal controls
of systems governed by impulsive integro-differential equations of mixed type. Wang\cite{wang2022backward} investigate the optimal control problems in terms of maximum principles and linear quadratic control problems of optimal control for forward stochastic Volterra integro-differential equations. 
 
 In this paper we focus on the following  Volterra  control system
\begin{align}\label{1.2}
\left\{\begin{array}{ll}
dX_t =b\Big(t,X_t ,\int_0^tk(t,s)X_sds,u_t,\int_0^tl(t,s)u_sds \Big)dt+\sigma\Big(t,X_t ,\int_0^tk(t,s)X_sds,u_t,\int_0^tl(t,s)u_sds  \Big)dW_t,\qquad 0\le t\le T,
\\X_0=x,
\end{array}\right.
\end{align}
to minimize the cost function
\begin{align*}
    J(u)=E\left[\int_0^T f\left(t,X_t,\int_0^tk(t,s)X_sds,u_t,\int_0^tl(t,s)u_sds \right)dt+g(X_T)\right].
\end{align*}

We call the corresponding stochastic differential equation (SDE in short)  as a generalized stochastic Volterra integral differential equation (SVIDE in short), which is a specific type of 
integral differential equations. This type of equation is influenced by the past information of both state process and control process from beginning to present.  

 It should be pointed out that most SVIDEs can not be written as stochastic Volterra integral equations with the following form
 \begin{align*}
X(t)=\varphi(t)+\int_0^t b(t, s, X(t),X(s)) d s+\int_0^t \sigma(t, s, X(t),X(s)) d W(s), \quad t \in[0, T] .
 \end{align*}
 For example, consider the following SDE, 
\begin{align*}
    dX_t=\left( X_t+{\rm sin}\left(\int_0^t X_s ds\right)\right)dt+\sigma(t)dW_t,
\end{align*}
where the drift term is nonlinear on the integral of state process.

We study the uniqueness of the solution of the SVIDE. Different from classical condition, to deal with the 
 term contain the past information, we use Gronwall inequality to $\sup_{0\le r\le t} E|\tilde{X}_r-X_r|^2$, which is a upper bound for $CE\left[\left|\frac{1}{t}\int_0^t\tilde{X}_sds-\frac{1}{t}\int_0^tX_sds\right|\right]^2$. Then a new stochastic maximum principle for control system (\ref{1.2}) is established. We define the Hamiltonian function 
and the adjoint equation to obtain the optimal system.
To study the properties of the adjoint equation, we prove the existence and uniqueness of the solution of the following equation
\begin{align*}
\left\{\begin{array}{ll}
-dy_t=h\left(t,y_t,z_t,E^{\mathcal{F}_t}\left[\int_t^Tk(s,t)a_1(s)y_sds\right],E^{\mathcal{F}_t}\left[\int_t^Tk(s,t)a_2(s)z_sds\right]\right)dt-z_tdW_t,
 \\
\\y_T=\xi ,
\end{array}\right.
\end{align*} which is singular at point 0. Compared with classical type investigated by El Karoui and Peng\cite{el1997backward}, we construct a contraction mapping under a new $\beta$-norm 
\begin{align*}
\|(Y,Z)\|_{\beta}=\sup_{0\le s\le T}Ee^{\beta s}| Y_s|^2+ E\int_0^Te^{\beta s} Z_s^2ds.
\end{align*}
Furthermore, we get the  necessary condition that the optimal control process should satisfy.
Consider the linear quadratic
case, which can be applied to a Volterra linear quadratic state regulator, we obtain the unique optimal control process for linear quadratic Volterra control system.

The rest of this paper is organized as follows. In section 2, we introduce a type of generalized SVIDE and prove the existence and the uniqueness of the solution of this type of equation. 
In section 3, we prove the stochastic maximum principle by introducing a kind of anticipated backward stochastic differential equations, and the existence and uniquenes of this kind of equations is proved.  
In the section 4, the linear quadratic case is investigated to illustrate the main result.

\section{A generalized  stochastic Volterra integral differential equation}

\qquad Let $(\Omega,\mathcal{F},\mathbb{P})$ be a probability space.    $\mathcal{F}_0\subset \mathcal{F}$ be a sub $\sigma$-algebra, and $\mathbb{F}=(\mathcal{F}_t)_{0\le t\le T}$ be the filtration generated by $\mathcal{F}_0$ and a $m$-dimensional standard Brownian motion $\textbf{W}=(W_t)_{0\le t\le T}$. We consider the following stochastic differential equation.

\begin{align}\label{masde}
\left\{\begin{array}{ll}
dX_t =b(t,X_t ,Y_t )dt+\sigma(t,X_t ,Y_t )dW_t,\qquad 0\le t\le T,
\\X_0=x,
\end{array}\right.
\end{align}
where $E|X_0|^2< \infty$, $b$ and $\sigma $ be measurable functions on $[0,T]\times \mathbf{R}^d\times \mathbf{R}^d$ with values in $\mathbf{R}^d$ and $\mathbf{R}^{d\times m}$, respectively. Here

\begin{align}\label{2.2}
Y_t=\left\{\begin{array}{ll}
\int_0^tk(t,s) X_s ds,&\textrm{$t>0$},
\\X_0,&\textrm{$t=0$},
\end{array}\right.
\end{align}
where $k(t,s)$ satisfies $\sup_{0\le t\le T}\int_0^t|k(t,s)|\le M$ for some constants $M>0$. It is obviously that $Y_t$ is continuous. This class of equations provides a description of the effect of past situations on the current situation. Assume that 
\begin{equation}\label{2.3}
    |b(t,x,y)|^2\lor|\sigma (t,x,y)|^2 \le L(1+|x|^2+|y|^2),\quad x,y\in \mathbf{R}^n,t\in [0,T],
\end{equation}
and
\begin{align}\label{2.4}
|b(t,x_1,y_1)-b(t,x_2,y_2)|^2&\lor|\sigma(t,x_1,y_1)-\sigma(t,x_2,y_2)|^2
\le L(|x_1-x_2|^2+|y_1-y_2|^2),\notag\\ &x_1,x_2,y_1,y_2\in \mathbf{R}^n,\quad t\in [0,T]
\end{align}
for some constant $L>0$ (where $|\sigma|^2=\sum |\sigma_{ij}|^2$).

\qquad Now we show the existence and uniqueness of the solution of equation (\ref{masde}).
~\\

\textbf{Lemma 2.1}\qquad If condition (\ref{2.3}) and (\ref{2.4}) holds, there exist a unique solution to equation (\ref{masde}).
~\\

\textbf{\emph{Proof}:}
Uniqueness. Let $X_t$ and $\tilde{X}_t$ be two solutions of the equation (\ref{masde}), $Y_t$ and $\tilde{Y}_t$ are corresponding moving average processes, and $X_0=Y_0=\tilde{X}_0=\tilde{Y}_0=x$.
Thus, we have

\begin{align*}
E|\tilde{X}_t-X_t|^2&=E\Big[\int_0^t b(s,\tilde{X}_s,\tilde{Y}_s)-b(s,X_s,Y_s)ds+\int_0^t \sigma(s,\tilde{X}_s,\tilde{Y}_s)-\sigma(s,X_s,Y_s)dW_s \Big]^2\notag\\
&\le2(T+1)LE\int_0^t |\tilde{X}_s-X_s|^2+|\tilde{Y}_s-Y_s|^2 ds\notag\\
&\le(2M^2+1)(T+1)L\int_0^t\sup_{0\le r\le s}E|\tilde{X}_r-X_r|^2ds
.\end{align*}
The last inequality holding is because
\begin{align*}
    E|\tilde{Y}_s-Y_s|^2=E\left|\int_0^s k(s,r)(\tilde{X}_r-X_r)dr\right|^2\le E\left\{\left[\int_0^s|k(s,r)|dr\right]\left[\int_0^s|k(s,r)||\tilde{X}_r-X_r|^2dr\right]\right\}\le M^2\sup_{0\le r\le s}E|\tilde{X}_r-X_r|^2.
\end{align*}
For every $\varepsilon>0$, there exits $\xi_t\in[0,t]$, such that 
\begin{align*}
    E|\tilde{X}_{\xi_t}-X_{\xi_t}|^2\ge\sup_{0\le r\le t} E|\tilde{X}_r-X_r|^2-\varepsilon,
\end{align*}
so that
\begin{align*}
\sup_{0\le r\le t} E|\tilde{X}_r-X_r|^2\le&
    E|\tilde{X}_{\xi_t}-X_{\xi_t}|^2+\varepsilon\\
    \le&(2M^2+1)(T+1)L\int_0^{\xi_t}\sup_{0\le r\le s}E|\tilde{X}_r-X_r|^2ds+\varepsilon\\
    \le&(2M^2+1)(T+1)L\int_0^t\sup_{0\le r\le s}E|\tilde{X}_r-X_r|^2ds+\varepsilon
.\end{align*}
Through the Gronwall's inequality and the arbitrariness of $\varepsilon$, we get $\sup_{0\le t\le T} E|\tilde{X}_t-X_t|^2=0$. Thus, the solution $X_t$ is unique.
~\\

Existence. Let
\begin{align*}
\left\{\begin{array}{ll}
X_t^{(k+1)}&=x+\int_0^t b(s,X_s^{(k)},Y_s^{(k)})dt+\int_0^t\sigma(s,X_s^{(k)},Y_s^{(k)})dW_s ,
\\\quad Y_t^{(k)}&=\int_0^t k(t,s)X_s^{(k)} ds,
\\\quad X_t^{(0)}&=x,
\end{array}\right.
\end{align*}
and
\begin{align*}
    u_t^{(k)}=\sup_{0\le r\le t}E\Big|X_r^{(k+1)}-X_r^{(k)}\Big|^2
.\end{align*}
 Similar to the proof of classical case, we get

\begin{align*}
u_t^{(k)}\le\frac{A^{k+1}t^{k+1}}{(k+1)!}
\end{align*}
for some constants $A>0$. Let $\lambda$ be Lebesgue measure on $[0,T]$, $0\le n < m$ and $m,n\to \infty$. Then  we have
~\\

\begin{align*}
\left\|X_{t}^{(m)}-X_{t}^{(n)}\right\|_{L^2(\lambda\times P)}\le\sum_{k=n}^{m-1}\left(\frac{A^{k+2} T^{k+2}}{(k+2) !}\right)^{\frac{1}{2}} \rightarrow 0 
.\end{align*}
Therefore, $\{X_t^{(n)}\}_{n\ge 0}$ is a Cauchy sequence in ${L^2(\lambda\times P)}$. Define
\begin{align*}
    X_t:=\lim_{n\to\infty} X_t^{(n)}
,\qquad
Y_t:=\lim_{n\to\infty} Y_t^{(n)}=\lim_{n\to\infty} \int_0^t k(t,s)X_s^{(n)}ds.
\end{align*}
Then $X_t$ and $Y_t$ are $\mathcal{F}_t$-measurable for all $t$. Since this holds for each $X_t^{(n)}$ and $Y_t^{(n)}$, thus $X_t$ is the solution of (\ref{masde}).

\section{The Maximum Principle}

\quad Consider the following control problem. The state equation is
\begin{align}
\left\{\begin{array}{ll}
dX_t =b\Big(t,X_t ,\int_0^tk(t,s)X_sds,u_t,\int_0^tl(t,s)u_sds \Big)dt+\sigma\Big(t,X_t ,\int_0^tk(t,s)X_sds,u_t,\int_0^tl(t,s)u_sds  \Big)dW_t,\qquad 0\le t\le T,
\\X_0=x,
\end{array}\right.
\end{align}
with the cost function
\begin{align}
    J(u)=E\left[\int_0^T f\left(t,X_t,\int_0^tk(t,s)X_sds,u_t,\int_0^tl(t,s)u_sds \right)dt+g(X_T)\right],
\end{align}
where $b(t,x,y,u,v)$ and $\sigma(t,x,y,u,v) $ are measurable functions on $\mathbf{R}\times \mathbf{R}^d\times \mathbf{R}^d\times \mathbf{R}^k\times \mathbf{R}^k$ with values in $\mathbf{R}^d$ and $\mathbf{R}^{d\times m}$, respectively.
$f(t,x,y,u,v)$ and $g(x) $ be measurable functions on $\mathbf{R}\times \mathbf{R}^d\times \mathbf{R}^d\times \mathbf{R}^k\times \mathbf{R}^k$ and $\mathbf{R}^d$, respectively, with values in $\mathbf{R}$. 
We denote by $\mathbb{U}$ the set of progressively measurable process
\textbf{u}$=(u_t)_{0\le t\le T}$ taking values in a given closed-convex set $\textbf{U}\subset \mathbb{R}^k$ and satisfying $E\int_0^T |u_t|^2 dt <\infty$.

To simplify the notation without losing the generality, we just consider the case $d=m=k=1$. We assume $u_t^{*}$ is the optimal control process, i.e.,
\begin{align*}
    J(u_t^*)=\min_{u_t\in \mathbb{U}}J(u_t)
.\end{align*}
For all $0<\varepsilon<1$, let 
\begin{align*}
    u_t^\varepsilon =(1-\varepsilon)u_t^*+\varepsilon\alpha_t\triangleq 
u^*_t+\varepsilon \beta_t,
\end{align*}
where $\alpha_t$ is any other admissible control.

We define the Hamiltonian function $H$ by
\begin{align}
    H(t,x,y,u,v,p,q)=b(t,x,y,u,v)p+\sigma (t,x,y,u,v)q+f(t,x,y,u,v)
.\end{align}
Denote
\begin{align*}
\phi^*(t)=\phi\left(t,X_t^*,\int_0^tk(t,s)X_s^*ds,u_t^*,\int_0^tl(t,s)u_s^*ds\right),
\end{align*}
for $\phi=b,\sigma,f,b_x,\sigma_x,f_x,b_y,\sigma_y,f_y,b_u,\sigma_u,f_u,b_v,\sigma_v,f_v$.
~\\

\textbf{Theorem 3.1} \qquad If $(u_t^*)_{0\le t\le T}$ is the optimal control process, $(X_t^*)_{0\le t\le T}$ and $(p_t,q_t)$ is the process satisfying \begin{align}\label{absde}
\left\{\begin{array}{ll}
-dp_t=\Big[b_x^*(t)p_t+E^{\mathcal{F}_t}\left[\int_t^Tk(s,t)b_y^*(s)p_sds\right]+\sigma_x^*(t)q_t+E^{\mathcal{F}_t}\left[\int_t^Tk(s,t)\sigma_y^*(s)q_sds\right]\\
\qquad\qquad +f_x^*(t)+E^{\mathcal{F}_t}\left[\int_t^Tk(s,t)f_y^*(s)ds\right]\Big]dt-q_tdW_t,
 \\
\\p_T=g_x(X_T^*) .
\end{array}\right.
\end{align}
Then we have
\begin{align}\label{3.25}
    \Big[H^*_u(t)+E^{\mathcal{F}_t}[\int_t^Tl(s,t)H^*_v(s)ds]\Big]\cdot(\alpha_t-u^*_t)\ge 0 ,
    \qquad \forall \alpha_t \in \mathbb{U},\quad d\lambda\otimes dP\quad a.s.
\end{align}
 for any  control process $\alpha_t$, where 
\begin{align*}
    H^*(t)=H\Big(t,X_t^* ,\int_0^tk(t,s)X_s^*ds,u_t^* ,\int_0^tl(t,s)u_s^*ds,p_t,q_t\Big).
\end{align*}
~\\

\textbf{Remark 3.2} \quad To investigate the adjoint equation (\ref{absde}), we consider a more general type of anticipated backward stochastic differential
equations:
\begin{align}\label{aabsde}
\left\{\begin{array}{ll}
-dy_t=h\left(t,y_t,z_t,E^{\mathcal{F}_t}\left[\int_t^Tk(s,t)a_1(s)y_sds\right],E^{\mathcal{F}_t}\left[\int_t^Tk(s,t)a_2(s)z_sds\right]\right)dt-z_tdW_t,
 \\
\\y_T=\xi .
\end{array}\right.
\end{align}
Without losing 
the generality, we assume $\sup_{0\le t\le T}\int_0^t|k(t,s)|ds=\sup_{0\le t\le T}\int_0^t|l(t,s)|ds=1$.  This type of anticipated backward stochastic differential equations may be singular even we assume $h(t,y,z,\tilde{y},\tilde{z})$ is Lipscitz continous, such as $k(t,s)=\frac{1}{t}$. We would deal with it in the following Lemma.

Anticipated backward stochastic differential equations was first studied by Peng and Yang \cite{peng2009anticipated}, they study the following type of equation:
\begin{equation*}
\begin{cases}-d Y_t=f\left(t, Y_t, Z_t, Y_{t+\delta(t)}, Z_{t+\zeta(t)}\right) d t-Z_t d W_t, & t \in[0, T] , \\ Y_t=\xi_t, & t \in[T, T+K], \\ Z_t=\eta_t, & t \in[T, T+K],\end{cases}
\end{equation*}
the unique solutions, a comparison theorem, and a duality
between them and stochastic differential delay equations are introduced. More properties of generalized anticipated
backward stochastic differential equations refer to Yang and Elliott  \cite{yang2013some}. 

~\\

\textbf{Lemma 3.3}\quad The anticipated backward stochastic differential equation (\ref{aabsde}) has the unique solution pair if the following conditions hold:
\begin{align*}
|h(t,y,z,\tilde{y},\tilde{z})|&\le M_1(|y|+|z|+|\tilde{y}|+|\tilde{z}|),\\
|h(t,y_1,z_1,\tilde{y_1},\tilde{z_1})-h(t,y_2,z_2,\tilde{y_2},\tilde{z_2})|&\le M_1(|y_1-y_2|+|z_1-z_2|+|\tilde{y}_1-\tilde{y}_2|+|\tilde{z_1}-\tilde{z_2}|),
\end{align*}
and
\begin{align}
a_1(s)\vee a_2(s)\le M_2, \qquad \forall s\in [0,T],\quad a.s.
\end{align}
for some constants $M_1, M_2>0$ satisfy $M_1M_2<8^{-1}T^{-\frac{1}{2}}$.

~\\

\textbf{\emph{Proof}:}\quad
Denote $\mathbb{H}_T^2\left(\mathbb{R}^d\right)$ is the space of all predictable processes $\phi: \Omega \times[0, T] \mapsto \mathbb{R}^d$ such that $\|\varphi\|^2=$ $\mathbb{E} \int_0^T\left|\varphi_t\right|^2 d t<+\infty.$
We define $\beta$-norm: $\|(Y,Z)\|_{\beta}=\sup_{0\le s\le T}Ee^{\beta s}| Y_s|^2+ E\int_0^Te^{\beta s} Z_s^2ds$ on $\mathbb{H}_T^2\left(\mathbb{R}^d\right)\times \mathbb{H}_T^2\left(\mathbb{R}^{d\times m}\right) $.
 For any $\mathcal{F}_t$-adapted 
 continuous process pair $(y_t^1,z_t^1), (y_t^2,z_t^2)$ with bounded $\beta$-norm, let
\begin{align*}
\left\{\begin{array}{ll}
-dY^i_t&=h\left(t,y^i_t,z^i_t,E^{\mathcal{F}_t}\left[\int_t^Tk(s,t)a_1(s)y^i_sds\right],E^{\mathcal{F}_t}\left[\int_t^Tk(s,t)a_2(s)z^i_sds\right]\right)dt-Z^i_tdW_t,
\\\quad Y^i_T&=\xi,
\end{array}\right.
\end{align*}
for $i=1,2$.

We denote
\begin{align*}
\delta \phi_t=\phi^1_t-\phi^2_t,
\end{align*}
for $\phi=Y,Z,y,z$, and
\begin{align*}
\delta h_t=&h\left(t,y^1_t,z^1_t,E^{\mathcal{F}_t}\left[\int_t^Tk(s,t)a_1(s)y^1_sds\right],E^{\mathcal{F}_t}\left[\int_t^Tk(s,t)a_2(s)z^1_sds\right]\right)\\&-h\left(t,y^2_t,z^2_t,E^{\mathcal{F}_t}\left[\int_t^Tk(s,t)a_1(s)y^2_sds\right],E^{\mathcal{F}_t}\left[\int_t^Tk(s,t)a_2(s)z^2_sds\right]\right).
\end{align*}
We can directly get
\begin{align*}
|\delta h_t|\le M_1\left(|\delta y_t|+|\delta z_t|\right)+M_1M_2\left(E^{\mathcal{F}_t}\left[\int_t^Tk(s,t)|\delta y_s|ds\right]+E^{\mathcal{F}_t}\left[\int_t^Tk(s,t)|\delta z_s|ds\right]\right).
\end{align*}
By it$\hat{\rm o}$'s formula,
\begin{align*}
d\left(e^{\beta t}\delta Y_t^2\right)=&\beta e^{\beta t}\delta Y_t^2dt+2e^{\beta t}\delta Y_td\delta Y_t+e^{\beta t}\left(d\delta Y_t\right)^2\notag\\
=&e^{\beta t}\left[\beta\delta Y_t^2-2\delta Y_t\delta h_t+\delta Z_t^2\right]dt+e^{\beta t}\delta Y_t\delta Z_tdW_t.
\end{align*}
Taking the expectation and integral, we have
\begin{align}\label{2222}
Ee^{\beta t}\delta Y_t^2+E\int_t^Te^{\beta s}\delta Z_s^2ds=&E\int_t^T e^{\beta s}\left[-\beta \delta Y_s^2+2\delta Y_s\delta h_s\right]ds\notag\\
\le&E\int_t^T e^{\beta s}\left[-\beta \delta Y_s^2+(c_1^{-1}+c_2^{-1})M_1|\delta Y_s|^2+c_1M_1|\delta y_s|^2+c_2M_1|\delta z_s|^2\right]ds\notag\\
&+2M_1M_2E\int_t^T e^{\beta s}|\delta Y_s|\left(E^{\mathcal{F}_s}\left[\int_s^Tk(r,s)|\delta y_r|dr\right]+E^{\mathcal{F}_s}\left[\int_s^Tk(r,s)|\delta z_r|dr\right]\right)ds.
\end{align}
Notice that 
\begin{align*}
2E\left[\int_t^Te^{\beta s}|\delta Y_s|E^{\mathcal{F}_s}[\int_s^T k(r,s)|\delta y_r|dr ]ds \right]=&2E\int_t^T\int_s^Tk(r,s)e^{\beta s}|\delta Y_s||\delta y_r|drds\\
\le& c_3^{-1}E\int_t^T\int_s^Tk(r,s)e^{\beta s}|\delta Y_s|^2drds+c_3E\int_t^T\int_s^Tk(r,s)e^{\beta s}|\delta y_r|^2drds ,
\end{align*}
and
\begin{align*}
\int_t^T\int_s^Tk(r,s)drds=\int_t^T\int_t^rk(r,s)dsdr=\le T-t\le T,
\end{align*}
so we have
\begin{align}\label{333333}
2E\left[\int_t^Te^{\beta s}|\delta Y_s|E^{\mathcal{F}_s}[\int_s^T k(r,s)|\delta y_r|dr ]ds \right]\le&c_3^{-1}T\sup_{t\le s\le T}Ee^{\beta s}|\delta Y_s|^2+c_3E\int_t^Te^{\beta r}|\delta y_r|^2dr\notag\\
\le &c_3^{-1}T\sup_{t\le s\le T}Ee^{\beta s}|\delta Y_s|^2+c_3T\sup_{t\le s\le T}Ee^{\beta s}|\delta y_s|^2.
\end{align}
In the same way
\begin{align}\label{44444}
2E\left[\int_t^Te^{\beta s}|\delta Y_s|E^{\mathcal{F}_s}[\int_s^T k(r,s)|\delta z_r|dr ]ds \right]\le c_3^{-1}T\sup_{t\le s\le T}Ee^{\beta s}|\delta Y_s|^2+c_3E\int_t^Te^{\beta r}|\delta z_r|^2dr .
\end{align}
Substitute (\ref{333333}) and (\ref{44444}) into (\ref{2222}), choose $\beta >(c_1^{-1}+c_2^{-1})M_1$ and let $c_1,c_2\to \infty$, we get
\begin{align*}
Ee^{\beta t}|\delta Y_t|^2+ E\int_t^Te^{\beta s}\delta Z_s^2ds\le 2M_1M_2c_3^{-1}\sup_{t\le s\le T}Ee^{\beta s}|\delta Y_s|^2+2c_3M_1M_2T\left[\sup_{t\le s\le T}Ee^{\beta s}|\delta y_s|^2+E\int_t^Te^{\beta s}|\delta z_s|^2ds\right].
\end{align*}
Then it is not difficult to get
\begin{align*}
(1-4M_1M_2c_3^{-1})\sup_{0\le s\le T}Ee^{\beta s}|\delta Y_s|^2+ E\int_0^Te^{\beta s}\delta Z_s^2ds\le 4c_3M_1M_2T\left[\sup_{0\le s\le T}Ee^{\beta s}|\delta y_s|^2+E\int_0^Te^{\beta s}|\delta z_s|^2ds\right],
\end{align*}
which shows
\begin{align*}
(1-4M_1M_2c_3^{-1})\|(\delta Y,\delta Z)\|_{\beta}\le 4c_3M_1M_2T\|(\delta y,\delta z)\|_{\beta}.
\end{align*}
Under the assumption $M_1M_2<8^{-1}T^{-\frac{1}{2}}$ and taking $c_3=T^{-\frac{1}{2}}$, we get the contraction mapping $\Phi:(y,z)\to(Y,Z)$ from $\mathbb{H}_{T,\beta}^2\left(\mathbb{R}^d\right)\times\mathbb{H}_{T,\beta}^2\left(\mathbb{R}^{d\times m}\right)$ 
onto itself and that there exists a fixed point, which is the unique continuous solution of
the  antici-
pated backward stochastic differential equation.

This completes the proof of Lemma 3.3.

~\\

\textbf{Lemma 3.4} \qquad Let $(u_t^*)_{0\le t\le T}$ is the optimal control process and $(X_t^*)_{0\le t\le T}$ be the corresponding state process, and $(p_t,q_t)$ is the adjoint process satisfying(\ref{absde}).
Then the $G\hat{a}teaux$ derivative of $J$ at $u^*_t$ in the direction $\beta_t$ is
\begin{align}
    \frac{d}{d\varepsilon}J(u_t^*+\varepsilon \beta_t)\Big|_{\varepsilon=0}=E\int_0^T \Big[H^*_u(t)+\int_t^Tl(s,t)H^*_v(s)ds\Big]\cdot \beta_t dt
.\end{align}

\textbf{\emph{Proof}:}\quad Let $X_t^*$ and $X_t^\varepsilon$ be the state process corresponding to $u_t^*$ and $u_t^\varepsilon$, respectively. Define $V_t$ by
\begin{align}
\left\{\begin{array}{ll}
dV_t=&\Big[b_x^*(t) V_t+b_y^*(t)\int_0^t k(t,s)V_s ds+b_u^*(t) \beta_t+b_v^*(t)\int_0^tl(t,s)\beta_sds\Big]dt\\
&+\Big[\sigma_x^*(t) V_t+\sigma_y^*(t)\int_0^tk(t,s) V_s ds+\sigma_u^*(t) \beta_t+\sigma_v^*(t)\int_0^tl(t,s)\beta_sds\Big]d W_t,
\\V_0=0.
\end{array}\right.
\end{align}
 It's easy to get
\begin{align*}
\sup_{0\le t\le T}\lim_{\varepsilon\to 0} E\Big[\frac{X_t^\varepsilon-X_t^*}{\varepsilon}-V_t\Big]^2=0.
\end{align*}
~\\
So we have that                              
\begin{align}\label{3.7}
\frac{J(u_t^\varepsilon)-J(u_t^*)}{\varepsilon}\to
E\left[\int_0^T\left(f_x^*(t)V_t+f_y^*(t)\int_0^t k(t,s)V_s ds+f_u^*(t)\beta_t+f_v^*(t)\int_0^tl(t,s)\beta_sds\right)dt+g_x(X_T^*)V_T\right]
,\end{align}
as $\varepsilon\to 0$.

 By It$\hat{\rm o}$'s formula, we have that
\begin{align}
d(p_tV_t)&=p_tdV_t+V_tdp_t+dp_tdV_t\notag \\
    &=p_t\Big[b_x^*(t) V_t+b_y^*(t)\int_0^tk(t,s) V_s ds+b_u^*(t) \beta_t+b_v^*(t)\int_0^tl(t,s)\beta_sds\Big]dt\notag \\
    &\quad -V_t\Big[b_x^*(t)p_t+E^{\mathcal{F}_t}[\int_t^Tk(s,t)b_y^*(s)p_sds]+\sigma_x^*(t)q_t  +E^{\mathcal{F}_t}[\int_t^Tk(s,t)\sigma_y^*(s)q_sds]+f_x^*(t)+E^{\mathcal{F}_t}[\int_t^Tk(s,t)f_y^*(s)ds]\Big]dt\notag \\
    &\quad+q_t\Big[\sigma_x^*(t) V_t+\sigma_y^*(t)\int_0^tk(t,s) V_s ds+\sigma_u^*(t) \beta_t+\sigma_v^*(t)\int_0^tl(t,s)\beta_sds\Big]dt+M_tdW_t \notag \\
    &=\Big[b_y^*(t)p_t\int_0^t k(t,s)V_s ds-V_tE^{\mathcal{F}_t}\int_t^Tk(s,t)b_y^*(s)p_sds
    +\sigma_y^*(t)q_t\int_0^t k(t,s)V_s ds-V_tE^{\mathcal{F}_t}[\int_t^Tk(s,t)\sigma_y^*(s)q_sds]\Big]d t\notag \\
    &\quad-\Big[f_x^*(t)V_t-V_tE^{\mathcal{F}_t}[\int_t^Tk(s,t)f_y^*(s)ds]\Big]dt+\Big[b_u^*(t)p_t+\sigma_u^*(t)q_t\Big]\beta_tdt\notag\\
    &\quad+\Big[b_v^*(t)p_t+\sigma_v^*(t)q_t\Big]\int_0^tl(t,s)\beta_sdsdt+M_t dW_t
,\end{align}
where $(M_t)_{0\le t\le T}$ is a $\mathcal{F}_t$ adapted process.

Consider
\begin{align}\label{3.14}
    Eg_x(X_T^*)V_T
    &=Ep_TV_T=E\int_0^T d(p_tV_t)+Ep_0V_0\notag\\
    &=E\int_0^T b_y^*(t)p_t\int_0^t k(t,s)V_s dsdt
     -E\int_0^T V_t\int_t^T\frac{1}{s}b_y^*(s)p_sdsdt\notag\\
    &\quad+E\int_0^T  \sigma_y^*(t)q_t\int_0^tk(t,s) V_s dsdt
     -E\int_0^T V_t\int_t^Tk(s,t)\sigma_y^*(s)q_sdsdt\notag\\
    &\quad -E \int_0^T f_x^*(t)V_tdt-E\int_0^T V_t\int_t^Tk(s,t)f_y^*(s)dsdt\notag\\
    &\quad +E\int_0^T \Big[b_u^*(t)p_t+\sigma_u^*q_t\Big]\beta_tdt+E\int_0^T\Big[b_v^*(t)p_t+\sigma_v^*(t)q_t\Big]\int_0^tl(t,s)\beta_sdsdt   
.\end{align}
By exchanging the order of integration, we have
\begin{align}
    \int_0^T b_y^*(t)p_t\int_0^t k(t,s)V_s dsdt=\int_0^T V_t\int_t^Tk(s,t)b_y^*(s)p_sdsdt,\label{3.15}\\
    \int_0^T  \sigma_y^*(t)q_t\int_0^tk(t,s) V_s dsdt
     =\int_0^T V_t\int_t^Tk(s,t)\sigma_y^*(s)q_sdsdt\label{3.16}
,\end{align}
and
\begin{align}\label{3.17}
    \int_0^T V_t\int_t^Tk(s,t)f_y^*(s)dsdt=\int_0^T f_y^*(t)\int_0^tk(t,s)V_sdsdt
.\end{align}
Substitute (\ref{3.15}), 
 (\ref{3.16}) into (\ref{3.14}), we get
\begin{align}\label{3.18}
  Eg_x(X_T^*)V_T
    &=-E \int_0^T f_x^*(t)V_tdt-E\int_0^T V_t\int_t^Tk(s,t)f_y^*(s)dsdt\notag\\
    &\quad +E\int_0^T \Big[b_u^*(t)p_t+\sigma_u^*(t)q_t\Big]\beta_tdt
+E\int_0^T\Big[b_v^*(t)p_t+\sigma_v^*(t)q_t\Big]\int_0^tl(t,s)\beta_sdsdt   .\end{align}
Then, substitute (\ref{3.18}) into (\ref{3.7}) and by (\ref{3.17}), we have
\begin{align}
    \frac{d}{d\varepsilon}J(u_t^*+\varepsilon \beta_t)\Big|_{\varepsilon=0}&=E\int_0^T\Big(f_x^*V_t+f_y^*int_0^t k(t,s)V_s ds+f_u^*\beta_t+f_v^*(t)\int_0^tl(t,s)\beta_sds\Big)dt\notag\\
    &\quad-E \int_0^T f_x^*(t)V_tdt-E\int_0^T V_t\int_t^Tk(s,t)f_y^*(s)dsdt\notag\\
    &\quad +E\int_0^T \Big[b_u^*(t)p_t+\sigma_u^*(t)q_t\Big]\beta_tdt+E\int_0^T\Big[b_v^*(t)p_t+\sigma_v^*(t)q_t\Big]\int_0^tl(t,s)\beta_sdsdt\notag\\
    &=E\int_0^T\Big[b_u^*(t)p_t+\sigma_t^*q_t+f_u^*(t)\Big]\beta_tdt+E\int_0^T\Big[b_v^*(t)p_t+\sigma_v^*(t)q_t+f_v^*(t)\Big]\int_0^tl(t,s)\beta_sdsdt\notag\\
    &=E\Big[\int_0^T H_u^*(t) \beta_t+H_v^*(t)\int_0^tl(t,s)\beta_sds\Big] dt\notag\\
    &=E\int_0^T \Big[H^*_u(t)+\int_t^Tl(s,t)H^*_v(s)ds\Big]\cdot \beta_t dt.
\end{align}
The last equality holding is because that 
 $$\int_0^TH_v^*(t)\int_0^tl(t,s)\beta_sdsdt=\int_0^T\beta_t\int_t^Tl(s,t)H_v^*(s)dsdt.$$ This completes the proof of Lemma 3.4.
~\\

Since $(u_t^*)_{0\le t\le T}$ is optimal control process, we have the inequality
\begin{align*}
    \frac{d}{d\varepsilon}J\Big(u_t^*+\varepsilon (\alpha_t-u_t^*)\Big)\Big|_{\varepsilon=0}\ge 0
.\end{align*}
By Lemma 3.4, we get
\begin{align*}
   E\int_0^T \Big[H^*_u(t)+\int_t^Tl(s,t)H^*_v(s)ds\Big]\cdot (\alpha_t-u^*_t) dt\ge 0
.\end{align*}
So
\begin{align*}
    E\Big[\mathbf{1}_A\big[H^*_u(t)+\int_t^Tl(s,t)H^*_v(s)ds\big]\Big]\cdot (\alpha_t-u^*_t)\ge 0,\quad \forall t\in[0,T],\quad \forall A\subset \mathcal{F}_t.
\end{align*}
To ensure adaptability, we can rewrite the above equation as
\begin{align*}
    E\Big[\mathbf{1}_A\big[H^*_u(t)+E^{\mathcal{F}_t}[\int_t^Tl(s,t)H^*_v(s)ds]\big]\Big]\cdot (\alpha_t-u^*_t)\ge 0,\quad \forall t\in[0,T],\quad \forall A\subset \mathcal{F}_t,
\end{align*}
and obtain that
\begin{align*}
    \Big[H^*_u(t)+E^{\mathcal{F}_t}[\int_t^Tl(s,t)H^*_v(s)ds]\Big]\cdot (\alpha_t-u^*_t) \ge 0,\qquad \forall t\in[0,T].
\end{align*}
This completes the proof of Theorem 3.1.

~\\

\textbf{Remark 3.5}  \qquad If the optimal control process $(u_t^*)_{0\le t\le T}$ takes values in the interior of the $\mathbb{U}$ , then we can replace (\ref{3.25})  with the following condition
\begin{align*}
    H_u^*(t)+E^{\mathcal{F}_t} 
 \Big[\int_t^Tl(s,t)H_v^*(s)ds\Big]=0
.\end{align*}
~\\

Thus, we give the optimal system
 \begin{align}
    \left\{\begin{array}{ll}
dX_t^*=H_p^*(t)dt+H_q^*(t)dW_t,\\\\
  -dp_t=\Big[H^*_x(t)+E^{\mathcal{F}_t}[\int_t^Tk(s,t)H^*_y(s)ds]\Big]dt -q_tdW_t,\\\\
  X_0^*=x,\quad  p_T=g_x(X_T^*),\\\\
  H_u^*(t)+E^{\mathcal{F}_t}
 [\int_t^Tl(s,t)H_v^*(s)ds]=0,
\end{array}\right.
\end{align}

where
\begin{align*}
    H^*(t)&=H\Big(t,X_t^* ,\int_0^tk(t,s)X_s^*ds,u_t^* ,\int_0^tl(t,s)u_s^*ds,p_t,q_t\Big),\\
    H(t,x,y,u,v,p,q)&=b(t,x,y,u,v)p+\sigma(t,x,y,u,v)q+f(t,x,y,u,v).
\end{align*}

\section{Linear quadratic case}

\quad In this section, we consider a linear quadratic (LQ in short) case, which can describe a moving average linear quadratic regulator problem. For simplicity, let $Y_t$ be the moving average process defined as (\ref{2.2}) and $v_t=\int_0^tl(t,s)u_sds$. The state process is defined as follows
\begin{align}\label{4.1}
    dX_t=\Big(A_tX_t+B_tY_t+C_tu_t+P_tv_t\Big)dt+\Big(D_tX_t+F_tY_t+H_tu_t+N_tv_t\Big)dW_t ,
\end{align}
with the cost function
\begin{align}\label{4.2}
    J(u)=\frac{1}{2}E\Big[\int_0^T (Q_tX_t^2+S_tY_t^2+R_tu_t^2)dt+GX_T^2\Big].
\end{align}
Here $G>0$ and $Q_t,S_t,R_t$ are positive functions.
~\\

Using the conclusions of Section 3, we can get the adjoint equation
\begin{align}\label{4.3}
    \left\{\begin{array}{ll}
-dp_t=\Big[A_tp_t+E^{\mathcal{F}_t}[\int_t^Tk(s,t)B_sp_sds]+D_tq_t+E^{\mathcal{F}_t}[\int_t^Tk(s,t)F_sq_sds]+Q_tX^*_t+E^{\mathcal{F}_t}[\int_t^Tk(s,t)S_sY^*_sds]\Big]dt\\
\\
 \qquad \qquad -q_tdW_t,
 \\
\\p_T=GX_T^* ,
\end{array}\right.
\end{align}
 and the optimal control process $u_t^*$ should satisfy
 \begin{align*}
     C_tp_t+H_tq_t+R_tu^*_t+E^{\mathcal{F}_t}[\int_t^Tl(s,t)P_sp_sds]+E^{\mathcal{F}_t}[\int_t^Tl(s,t)N_sq_sds]=0,
 \end{align*}
i.e.,
\begin{align}\label{4.4}
    u_t^*=-R_t^{-1}\left(C_tp_t+H_tq_t+E^{\mathcal{F}_t}[\int_t^Tl(s,t)P_sp_sds]+E^{\mathcal{F}_t}[\int_t^Tl(s,t)N_sq_sds]\right).
\end{align}

\textbf{Theorem 4.1} The function $u_t^*=-R_t^{-1}\left(C_tp_t+H_tq_t+E^{\mathcal{F}_t}[\int_t^Tl(s,t)P_sp_sds]+E^{\mathcal{F}_t}[\int_t^Tl(s,t)N_sq_sds]\right),\quad t\in[0,T]$ is the unique optimal control for moving average LQ problem (\ref{4.1}), (\ref{4.2}), where $(p_t,q_t)$ is defined by equality (\ref{4.3}). 
~\\

\textbf{\emph{Proof}:} \quad We now prove $u_t^*$ is the optimal control. For any $\tilde{u}_t\subset \mathbb{U}$, let $(\tilde{X}_t,\tilde{Y}_t,\tilde{v}_t)$ and $(X_t^*,Y_t^*,v_t^*)$ are processes corresponding to $\tilde{u}_t$ and $u_t^*$, respectively. We have that 
\begin{align*}
    d(\tilde{X}_t-X_t^*)=&[A_t(\tilde{X}_t-X^*_t)+B_t(\tilde{Y}_t-Y^*_t)+C_t(\tilde{u}_t-u^*_t)+P_t(\tilde{v}_t-v^*_t)]dt\\
    &+[D_t(\tilde{X}_t-X^*_t)+F_t(\tilde{Y}_t-Y^*_t)+H_t(\tilde{u}_t-u^*_t)+N_t(\tilde{v}_t-v^*_t)]dW_t.
\end{align*}
Consider
\begin{align}\label{4.5}
    dp_t(\tilde{X}_t-X_t^*)=&p_td(\tilde{X}_t-X^*_t)+(\tilde{X}_t-X^*_t)dp_t+dp_td(\tilde{X}_t-X_t^*)\notag\\
    =&p_t\left[A_t(\tilde{X}_t-X^*_t)+B_t(\tilde{Y}_t-Y^*_t)+C_t(\tilde{u}_t-u^*_t)+P_t(\tilde{v}_t-v^*_t)\right]dt\notag\\
    &-(X_t-X_t^*)\Bigg[A_tp_t+E^{\mathcal{F}_t}\left[\int_t^Tk(s,t)B_sp_sds\right]+D_tq_t+E^{\mathcal{F}_t}\left[\int_t^Tk(s,t)F_sq_sds\right]\notag\\
    &\qquad\qquad\qquad+Q_tX^*_t+E^{\mathcal{F}_t}\left[\int_t^Tk(s,t)S_sY^*_sds\right]\Bigg]dt\notag\\  
    &+q_t\left[D_t(\tilde{X}_t-X^*_t)+F_t(\tilde{Y}_t-Y^*_t)+H_t(\tilde{u}_t-u^*_t)+N_t(\tilde{v}_t-v^*_t)\right]dt+M_tdW_t\notag\\
    =&p_tB_t(\tilde{Y}_t-Y_t^*)dt-(\tilde{X}_t-X_t^*)
E^{\mathcal{F}_t}\left[\int_t^Tk(s,t)B_sp_sds\right]dt\notag\\
    &+q_tF_t(\tilde{Y}_t-Y^*_t)dt-(\tilde{X}_t-X_t^*)
E^{\mathcal{F}_t}\left[\int_t^Tk(s,t)F_sq_sds\right]dt\notag\\
    &+(C_tp_t+H_tq_t)(\tilde{u}_t-u^*_t)dt-Q_tX_t^*(\tilde{X}_t-X_t^*)dt\notag\\
    &+p_tP_t(\tilde{v}_t-v^*_t)dt+q_tN_t(\tilde{v}_t-v^*_t)dt\\
    &-(\tilde{X}_t-X_t^*)E^{\mathcal{F}_t}\left[\int_t^Tk(s,t)S_sY^*_sds\right]dt+M_tdW_t,
\end{align}
where $(M_t)_{0\le t\le T}$ is a $\mathcal{F}_t$ adapted process. By exchanging the order of integration, we get
\begin{align}\label{4.6}
    \int_0^Tp_tB_t(\tilde{Y}_t-Y_t^*)dt=\int_0^Tp_tB_t\int_0^tk(t,s)(\tilde{X}_s-X_s^*)dsdt=\int_0^T(\tilde{X}_t-X_t^*)\int_t^Tk(s,t)_sp_sdsdt.
\end{align}
In the same way, we have
\begin{align}
    \int_0^Tq_tF_t(\tilde{Y}_t-Y_t^*)dt=\int_0^T(\tilde{X}_t-X_t^*)\int_t^Tk(s,t)F_sq_sdsdt,\label{4.7} \\ 
     \int_0^Tp_tP_t(\tilde{v}_t-v_t^*)dt=\int_0^T(\tilde{u}_t-u_t^*)\int_t^Tl(s,t)P_sp_sdsdt,\label{4.71} \\ 
    \int_0^Tq_tN_t(\tilde{v}_t-v_t^*)dt=\int_0^T(\tilde{u}_t-u_t^*)\int_t^Tl(s,t)N_sq_sdsdt,\label{4.72}
\end{align}
and
\begin{align}\label{4.8}
    \int_0^T(\tilde{X}_t-X_t^*)\int_t^Tk(s,t)S_sY^*_sdsdt=\int_0^TS_tY^*_t(\tilde{Y}_t-Y_t^*)dt.
\end{align}

Taking integral for the (\ref{4.5}) from $0$ to $T$ and taking the expectation, through (\ref{4.4}), (\ref{4.6}), (\ref{4.7}),(\ref{4.71}),(\ref{4.72}) and (\ref{4.8}) we have
\begin{align}
EGX_T^*(\tilde{X}_T-X_T^*)=&Ep_T(\tilde{X}_T-X_T^*)\notag\\
    =&E\int_0^Tdp_t(\tilde{X}_T-X_T^*)\notag\\
    =&-E\int_0^T\Big[R_tu_t^*(\tilde{u}_t-u^*_t)+Q_tX_t^*(\tilde{X}_t-X^*_t)+S_tY^*_t(\tilde{Y}_t-Y_t^*)\Big]dt.
\end{align}
So that
\begin{align*}
J(\tilde{u}_t)-J(u_t^*)=&\frac{1}{2}E\int_0^T \Big[Q_t(\tilde{X}_t^2-X_t^{*2})+S_t(\tilde{Y}_t^2-Y_t^{*2})+R_t(\tilde{u}_t^2-u_t^{*2})\Big]dt \\
    &+\frac{1}{2}EG(\tilde{X}_T^2-X_T^{*2}) \\
    =& \frac{1}{2}E\int_0^T \Big[Q_t(\tilde{X}_t^2-X_t^{*2})-2Q_tX_t^*(\tilde{X}_t-X^*_t)+S_t(\tilde{Y}_t^2-Y_t^{*2})\\
     &\qquad\qquad-2S_tY^*_t(\tilde{Y}_t-Y_t^*)
   +R_t(\tilde{u}_t^2-u_t^{*2})-2R_tu_t^*(\tilde{u}_t-u^*_t)\Big]dt \\
    \ge& 0.
\end{align*}
This shows that $u_t^*$ is an optimal control. 

Then we prove $u^*_t$ is unique, assume that both $u_t^{*,1}$ and $u_t^{*,2}$ are optimal controls, $X_t^1$ and $X_t^2$ are corresponding state processes, respectively. It is easy to get $\frac{X_t^1+X_t^2}{2}$ is the corresponding state process to $\frac{u_t^{*,1}+u_t^{*,2}}{2}$. We assume there exist constants $\delta>0, \alpha\ge 0$, such that $R_t\ge \delta$ and 
\begin{align*}
    J(u_t^{*,1})=J(u_t^{*,2})=\alpha.
\end{align*}
Using the fact $a^2+b^2=2[(\frac{a+b}{2})^2+(\frac{a-b}{2})^2]$, we have that

\begin{align*}
    2\alpha=&J(u_t^{*,1})+J(u_t^{*,2})\\
           =&\frac{1}{2}E\int_0^T \Big[Q_t(X_t^1X_t^1+X_t^2X_t^2)+S_t(Y_t^1Y_t^1+Y_t^2Y_t^2)+R_t(u_t^{*,1}u_t^{*,1}+u_t^{*,2}u_t^{*,2})\Big]dt\\
           &+\frac{1}{2}EG(X_T^1X_T^1+X_T^2X_T^2)\\
           \ge&E\int_0^T \Big[Q_t\Big(\frac{X_t^1+X_t^2}{2}\Big)^2+S_t\Big(\frac{Y_t^1+Y_t^2}{2}\Big)^2+R_t\Big(\frac{u_t^{*,1}+u_t^{*,2}}{2}\Big)^2\Big]dt\\
           &+EG\Big(\frac{X_T^1+X_T^2}{2}\Big)^2+E\int_0^TR_t\Big(\frac{u_t^{*,1}-u_t^{*,2}}{2}\Big)^2dt\\
           =&2J\Big(\frac{u_t^{*,1}+u_t^{*,2}}{2}\Big)+E\int_0^TR_t\Big(\frac{u_t^{*,1}-u_t^{*,2}}{2}\Big)^2dt\\
           \ge&2\alpha+\frac{\delta}{4}E\int_0^T|u_t^{*,1}-u_t^{*,2}|^2dt.
\end{align*}
Thus, we have 
\begin{align*}
    E\int_0^T|u_t^{*,1}-u_t^{*,2}|^2dt\le 0,
\end{align*}
which shows that $u_t^{*,1}=u_t^{*,2}$.

\section*{Acknowledgments}

\qquad The authors acknowledge the financial support from the National Science Foundation of China (grant no. 11871244).

\bibliography{main}

\end{document}